\documentclass[12pt]{amsart}

\usepackage{graphicx}
\usepackage{amsfonts}
\usepackage{epsf}
\usepackage{amssymb}
\usepackage{amsmath}
\usepackage{amscd}
\usepackage{hyperref}
\usepackage{color}
\usepackage{bbm}
\usepackage{caption}
\usepackage{subcaption}

\newtheorem{theorem}{Theorem}[section]
\newtheorem*{utheorem}{Main Theorem}

\theoremstyle{remark}
\newtheorem{definition}[theorem]{Definition}

\newtheorem{remark}{Remark}[section]

\newcommand{\yes}{\mathrm{Yes}}
\newcommand{\no}{\mathrm{No}}

\input xy
\xyoption{all}

\makeatletter
\@namedef{subjclassname@2020}{\textup{2020} Mathematics Subject Classification}
\makeatother

\begin{document}
\title[NP-completeness of the $\ell_1$-embedding problem for simple graphs as SIGs]{NP-completeness of the $\ell_1$-embedding problem for simple graphs as sphere-of-influence graphs}
\begin{abstract}
In graph theory an interesting question is whether for a fixed choice of $p\in [0,\infty]$, all simple graphs appear as sphere-of-influence graphs in some Euclidean space with respect to the $\ell_p$ metric. The answer is affirmative for $p=\infty$, negative for any $p\in (0,\infty)$, and unknown for $p=1$. The result of this work shows that for the case of $p=1$, this embeddability question is a (Promise) NP-Complete problem.  
\end{abstract}
\author{Stanislav Jabuka}
\email{jabuka@unr.edu}
\address{Department of Mathematics and Statistics, University of Nevada, Reno NV 89557, USA.}
\maketitle

\section{Introduction}
\subsection{Background} Let $(M,d)$ be a metric space and $X=\{X_1,\dots, X_n\}\subset M$ be a finite subset. The {\em radius of influence $r_i$} of $X_i$ is $r_i = \min_{j\ne i} d(X_i, X_j)$, and the {\em sphere of influence} of $X_i$ is $S_i = \{x\in M\, |\, d(X_i, x)<r_i\}$.
The {\em sphere-of-influence graph}, or SIG for short, of the subset $X\subset M$ and induced by the metric $d$, is the graph with vertices $X_1, \dots, X_n$, two of which share an edge if and only if $d(X_i, X_j)<r_i+r_j$ (or equivalently, if their spheres of influence overlap). A central question in the theory of SIGs is whether any simple graph $G$ is isometric to a SIG in $(\mathbb R^m,d_p)$ with the dimension $m\ge 1$ arbitrary, and with $d_p$ being the $\ell_p$-metric 
$$d_p(x,y) = \left\{ 
\begin{array}{cl}
\left(|x_1-y_1|^p+ \dots+|x_n-y_n|^p\right)^{1/p} &\quad ; \quad p\in [1,\infty), \cr 
\max \{|x_1-y_1|, \dots,|x_n-y_n|\}& \quad ; \quad p=\infty.  	
\end{array} 
\right. 
$$
The answer is `yes' if $p=\infty$ \cite{DezaLaurentBook, MichaelQuint1}, it is `no' if $p\in (1,\infty)$ \cite{JabukaMirbagheri} and is unknown if $p=1$. We show that the case of $p=1$ leads to a (Promise) NP-Complete problem, lending credence to the difficulty of the problem. 

\subsection{$\ell_1$-embeddability and the cut cone $CUT_n$} Consider a finite metric space $(V_n,d)$ with $V_n=\{1,\dots, n\}$. A {\em cut of $V_n$} is a subset $C\subset V_n$, which we view as ``cutting'' $V_n$ into the disjoint pieces $V_n = C\cup (V_n-C)$. A cut creates a notion of distance on $V_n$, with two elements that lie in the same subset induced by the cut, either $C$ or $V_n-C$, being distance $0$ apart, and with two elements lying in opposite components of the cut being distance $1$ apart. This simple idea defines a semi-metric $\delta_C$ on $V_n$ by 
\begin{equation} \label{Definition of the cut metric}
\delta_C(i,j) = \left\{ 
\begin{array}{cl} 
1 & \quad ; \quad |C\cap \{i,j\}| = 1, \cr 
0 & \quad ; \quad |C\cap \{i,j\}| = 0 \text{ or } 2.  	
\end{array} 
\right. 	
\end{equation} 
We call $\delta_C$ the {\em cut-metric} corresponding to $C$. Obviously $\delta_C = \delta_{V_n-C}$. Let $CUT_n$ denote the {\em cut cone}, defined as 
\begin{align*}
\text{CUT}_n & = \left\{ \sum_{C\subset V_n} w_C \, \delta_C \,\big|\, w_C\ge 0\right\}. 
\end{align*}  
A metric $d$ on $V_n$ lies in the cut cone if, by definition, it can be written as a linear combination of cut-metrics with non-negative coefficients:
$$d = \sum _{C\subset V_n} w_C\cdot \delta_C, \quad w_C\ge 0 \text{ for all } C\subset V_n.$$
The following beautiful theorem was proved in \cite{Assouad}.
\begin{theorem} \label{Theorem about the equivalence of ell1-embeddability and belonging to the cut cone}
The metric space $(V_n,d)$ is isometric to a subspace of $(\mathbb R^m,d_1)$ if and only if $d\in CUT_n$. 
\end{theorem}
This theorem shifts the question of embeddability of a finite metric space $(V_n,d)$ into $(\mathbb R^m,d_1)$, to the question of membership of $d$ in the cut cone $CUT_n$. The latter is known to be NP-Complete decision problem:
\begin{theorem} [Avis-Deza \cite{AvisDeza}] The problem of determining if a metric $d$ on $V_n$ belongs to the cut cone $CUT_n$, is an NP-Complete decision problem. 
\end{theorem}
Building on this, and using the notions of promise decision problems and polynomial reduction, we prove what is the main result of this paper. 
\begin{utheorem} The problem of determining if a simple graph $G$ with vertices $V_n = \{1,\dots, n\}$, admits an embedding $\varphi: V_n\to \mathbb R^m$ so that the sphere-of-influence graph of $\{\varphi(1), \dots, \varphi(n)\}$ with respect to the $\ell_1$-metric $d_1$ on $\mathbb R^m$ is isomorphic to $G$, is a Promise NP-Complete problem. 
\end{utheorem}
\begin{remark} 
The notion of a {\em Promise NP-Complete} problem is a technical term explained in Section \ref{Section on all that P vs NP stuff}. For now it suffices to say that every Promise NP-Complete problem reduces in polynomial time to any NP-Complete problem, and that therefore both types of problems are of equivalent computational complexity. 
\end{remark} 
\section{Background on Computational Complexity} \label{Section with background material}
The exposition in this section follows standard texts in computation complexity theory. For decision/language problems we benefited from \cite{AroraBoaz, Hemaspaandra, Papa, Sipser}, while for promise problems we followed \cite{Goldreich, GoldreichBook}.
\subsection{Promise and Decision Problems} 
As is customary in computational complexity theory, we let $\Sigma = \{0,1\}$ denote the ``universal alphabet'' and $\Sigma ^\ast =\{\varepsilon\}\cup \left( \cup _{n\ge 1}\Sigma^{\times n}\right)$ be its Kleene star, with $\varepsilon$ the empty word. 

A {\em promise problem $B$} is a pair of disjoint subsets $B_\yes, B_\no \subset \Sigma^\ast$, with the possibility that $B_\yes \cup B_\no \ne \Sigma^\ast$. The {\em domain of the promise} or simply the {\em promise}, denoted Dom$(B)$, is Dom$(B) = B_\yes \cup B_\no$. Promise problems were introduced in \cite{EvenSelmanYacobi} with a view toward applications in Public-Key Cryptography. They have since received considerable attention and applications elsewhere, see for instance \cite{GoldreichBook}. 

A {\em decision problem $A$} is a promise problem for which Dom$(A) = A_\yes \cup A_\no =\Sigma^\ast$. A decision problem is also referred to as a {\em language recognition problem} with the language $L$ being $L = A_\yes$, in which case the decision to make is whether word $x\in \Sigma^\ast$ belongs to the language $L$. 

We will assume that the reader is familiar with the concept of a {\em deterministic Turing machine}, see for example \cite[Chapter 2]{Papa}.  
\subsection{P versus NP, and Promise P versus Promise NP} \label{Section on all that P vs NP stuff}
Following \cite[Definitions 2.4 and 2.5]{GoldreichBook}, we make this definition.
\begin{definition} \label{P and NP for Decision Problems}
Let $A$ be a decision problem. 
\begin{itemize}
\item[(i)] We say that $A$ is of {\em polynomial complexity}, or in {\em  class P}, if there exists a deterministic Turing machine $M$ on $\Sigma^\ast$ that halts in polynomial time, and accepts an $x\in \Sigma^\ast$ if and only if $x\in A_\yes$. 
\item[(ii)] $A$ is {\em non-deterministically polynomial} or in {\em class NP} if there exists a polynomial $p$ and a polynomial-time algorithm $V$ such that for every $x\in A_\yes$ there exists a $y\in \Sigma^\ast$ of length at most $p(|x|)$ and with $V(x,y)=1$ (such a $y$ is called an NP-witness for $x$). Additionally, for every $x\in A_\no$ and for any $y\in \Sigma^\ast$, $V(x,y)=0$ (saying that no witness $y$ can accept an element in $A_\no$). 
\end{itemize}
\end{definition}
Similar notions of computational complexity can be formulated for promise problems, see \cite[Definition 2.31.]{GoldreichBook}.
\begin{definition} \label{P and NP for Promise Problems}
Let $B$ be a promise problem. 
\begin{itemize}
\item[(i)] We say that $B$ is of {\em promise polynomial complexity}, or in {\em  class Promise P}, if there exists a deterministic Turing machine $M$ on $\Sigma^\ast$ that halts in polynomial time for any $x\in \mathrm{Dom}(B)$ and accepts an $x\in \mathrm{Dom}(B)$ if and only if $x\in B_\yes$. There is no requirement for the behavior of $M$ on $\Sigma^\ast -\mathrm{Dom}(B)$.  
\item[(ii)] $B$ is {\em promise non-deterministically polynomial} or in {\em class Promise NP} if there exists a polynomial $p$ and a polynomial-time algorithm $V$ such that for every $x\in B_\yes$ there exists a $y\in \Sigma^\ast$ of length at most $p(|x|)$ and with $V(x,y)=1$. Additionally, for every $x\in A_\no$ and for any $y\in \Sigma^\ast$, $V(x,y)=0$. No requirements are made for how the algorithm $V$ behaves on $(x,y)$ with $x\notin \mathrm{Dom}(B)$.   
\end{itemize}
\end{definition}
\begin{remark} \label{Remark about NP-Complete and  Promise NP-Complete}
Decision problems are special cases of promise problems, and Definitions \ref{P and NP for Decision Problems} and \ref{P and NP for Promise Problems} are compatible in the sense that if a decision problem is viewed as a promise problem, then it belongs to the class Promise P if and only if it belong to the class P, and similarly for NP. 
\end{remark} 
\subsection{Karp Reductions} For our purposes the more specialized Karp reductions will suffice, in favor of simplicity, compared to the more general Cook reductions. Karp reductions were introduced by Karp \cite{Karp} in 1972 with the goal of reducing a large variety of computational problems to language recognition problems. 
Our definition expands \cite[2.11]{GoldreichBook} from decision to promise problems. 
\begin{definition}  \label{Definition of Karp reduction}
Consider a pair of promise problems 
$$A: \,A_\yes \cup A_\no \subset \Sigma ^\ast, \qquad \text{ and } \qquad B:\, B_\yes \cup B_\no \subset \Sigma ^\ast.$$
A  {\em promise-preserving Karp reduction} from $A$ to $B$ is a function $f:\Sigma^\ast \to \Sigma ^\ast$, computable in polynomial time, such that:
\begin{itemize}
\item[(i)] $f\left(\mathrm{Dom}(A) \right) \subset \mathrm{Dom}(B)$.
\item[(ii)] $x\in A_\yes$ if and only if $f(x) \in B_\yes$. 
\end{itemize} 
\end{definition}
If such a function exists, we write 
$$A\le_p B,$$
with the subscript ``p'' referring to ``polynomial time''.  Points (i) and (ii) in the preceding definition imply that $x\in A_\no$ if and only if $f(x)\in B_\no$.  The notation $A\le_p B$ is meant to suggest that problem $B$ is at least as hard as problem $A$, up to polynomial time. For the same deterministic Turing machine $M$ that accepts/rejects elements  from Dom$(B)$, can be used to accept/reject elements in Dom$(A)$ by applying $M$ to $f(x)$ for any choice of $x\in \mathrm{Dom}(A)$.  

Since decision problems are viewed here as special cases of promise problems, the notion of promise-preserving Karp reduction applies to them as well, with condition (i) from Definition \ref{Definition of Karp reduction} becomes trivially true. Since that condition is the promise-preserving condition, we simply refer to {\em Karp reduction} in the case of decision problems. 
\vskip1mm
A decision problem $B$ is called {\em NP-Hard} if $A\le _p B$ for any NP-Problem $A$, and it is called {\em NP-Complete} if it is NP and NP-Hard. The first ever NP-Complete problem was the SAT problem discovered by Stephen Cook \cite{Cook} in 1971 (and independently discovered by Leonid Levin from the Soviet Union in 1973). By now there are many dozens of known NP-Complete problems, see for instance \cite{GareyJohnson, GoldreichBook2}.  

Similarly, a promise problem $B$ is called  {\em Promise NP-Hard} if $A\le _p B$ for any Promise NP Problem $A$, and it is called {\em Promise NP-Complete} if it is Promise NP and Promise NP-Hard.

\begin{remark}
From the definitions alone, it is be obvious that if $A$ and $B$ are two decision problems of which $A$  Karp reduces to $B$ and $B$ is either of class  P or class NP, then so is $A$. Similarly, if $A$ and $B$ are two promise problems of which $A$ promise-preserving Karp reduces to $B$ and $B$ is either of class Promise P or class Promise NP, then so is $A$.  
\end{remark}
\section{Proof of the Main Theorem}
Fix $n\ge 3$, let $V_n=\{1,\dots,n\}$ and let $CUT_n$ be the cut cone of metrics on $V_n$. Consider the two problems below, both of which we view as promise problems, even as the first is also a decision problem: 
\begin{itemize}
\item[{\bf 1.}] {\bf Problem A: Cut-Cone Membership Decision Problem.} \\
Choose $\Sigma^\ast$ to be the set of all metrics $d$ on $V_n$, and consider the decision problem $A_\yes \cup A_\no =\Sigma^\ast$ with \\  
\hspace*{10mm} $A_\yes = \{d \, |\, d\in CUT_n\}$,\\
\hspace*{10mm} $A_\no = \{d \, |\, d\notin CUT_n\}$.\\
\item[{\bf 2.}] {\bf Problem B: Cut-Cone Membership with SIG-Promise Problem.}\\
Choose $\Sigma^\ast$ to be the set of all pairs $(G,d)$ with $G$ a simple graph with vertices $V_n$, and consider the promise problem $B_\yes \cup B_\no \subset \Sigma^\ast$ with \\  
\hspace*{10mm} $B_\yes = \{d \, |\, d \text{ is a SIG metric on }V_n \text{ for $G$, and } d\in CUT_n\}$,\\
\hspace*{10mm} $B_\no = \{d \, |\, d \text{ is a SIG metric on }V_n \text{ for $G$, and } d\notin CUT_n\}$.\\
\end{itemize}

Of these, the first problem is known to be NP-Complete \cite{AvisDeza}, which by Remark \ref{Remark about NP-Complete and  Promise NP-Complete} means it is also  Promise NP-Complete. We proceed to show that each can be Karp reduced to the other in polynomial time, in a promise-preserving way. 
\vskip3mm
\noindent {\bf A$\le_p$B }
Given $d\in A_\yes$, that is, given a metric $d$ on $V_n$ with $d\in CUT_n$, let $G_d$ be the sphere-of-influence graph associated to $d$. This graph can be computed in polynomial time in $n$ as it requires finding the $n$ radii of influence $r_i =\min _{j\ne i}d(i,j)$, $i=1,\dots, n$, and then computing the ${n\choose 2}$ quantities $r_i+r_j-d(i,j)$, $1\le i <j \le n$ (which if positive lead to an edge between $i$ and $j$, and otherwise don't). It is trivial to see that $d$ is a SIG metric for $G_d$ so that $(G_d,d)\in B_\yes$. 

Similarly, if $d\in A_\no$, then with the same $G_d$ as above, $d$ is a SIG metric on $V_n$ for $G_d$ but $d\notin CUT_n$, and so $(G_d,d)\in B_\no$. 

In summary, the function $d\mapsto (G_d,d)$ is a polynomial time, promise-preserving Karp reduction from Problem A to Problem B. 
\vskip3mm
\noindent{\bf B$\le_p$A } 
If $(G,d)\in B_\yes$ then $d\in CUT_n$ and therefore $d\in A_\yes$, while if $(G,d)\in B_\no$ then $d\notin CUT_n$ and thus $d\in A_\no$. The forgetful function $(G,d) \to d$ is a polynomial time, promise-preserving Karp reduction from Problem B to Problem A. This completes the proof of the Main Theorem. 
\bibliographystyle{plain}
\bibliography{bibliography}
\end{document}